\documentclass[10pt,twoside]{article}
\usepackage{a4}
\usepackage{amssymb,amsmath,amsthm,latexsym}
\usepackage{amsfonts}
\usepackage{graphicx}

\begin{document}

\label{'ubf'}
\setcounter{page}{1}

\markboth{\hspace*{-9mm}\centerline{\footnotesize \sc
   The number of distinct element orders in a finite group}}
   {\centerline{\footnotesize \sc Manoj Kumar Singh}\hspace*{-9mm}}

\vspace*{-2cm}

\begin{center}
{
{\Large \textbf{\sc The number of distinct element orders in a finite group}}
\\

\medskip

{\sc Manoj Kumar Singh}\\
{\footnotesize Department of Mathematics and Statistics, St.\ Xavier's College, Ranchi}\\
{\footnotesize affiliated to Ranchi University, Ranchi, Jharkhand 834001, India}\\
{\footnotesize ORCID: 0009-0002-4054-6537}\\
{\footnotesize e-mail: {\it manojkumarsingh@sxcran.org}}

\medskip

{\sc Sumant Kumar}\\
{\footnotesize Department of Mathematics, Gossner College, Ranchi, Jharkhand, India}\\
{\footnotesize ORCID: 0000-0001-7466-289X}\\
{\footnotesize e-mail: {\it sumant.math@gcran.org}}
}
\end{center}

\thispagestyle{empty}

\hrulefill

\begin{abstract}
{\footnotesize
For a finite group $G$ let $\omega(G)$ denote its spectrum, the set of orders of its elements, and set
$\eta(G)=|\omega(G)|$, the number of distinct element orders. The set $\omega(G)$ has been studied
intensively, but its cardinality has received little systematic attention as an invariant in its own
right. We develop the theory of $\eta$. After describing its behaviour under subgroups, quotients,
sections, direct products and Frobenius extensions, we prove the two-sided bound
$1+\sum_{p\mid|G|}v_p(\exp G)\le\eta(G)\le\tau(\exp G)$,
in which the lower equality characterises the CP-groups, those all of whose elements have prime-power
order, and the upper equality the groups whose spectrum realises every divisor of the exponent;
nilpotent groups lie at the upper end, giving $\eta(G)=\tau(\exp G)$. We also prove
$\eta(G)\le\tau(|G|)$ with equality only for cyclic groups, and $\eta(G)\ge\pi(|G|)+1$ with equality
only when every nontrivial element has prime order. We show that $\eta(G)\le3$ forces solvability by an
argument resting only on Burnside's theorem, and, combining this with the recognition of $A_5$ by its
spectrum, that $A_5$ is the unique non-solvable group with $\eta(G)=4$. For the symmetric and
alternating groups we record self-contained proofs of closed criteria for membership in the spectrum
through a single arithmetic function, placing the associated counting sequences in the present
framework, and we determine the extreme values of $\eta$ over all groups of a fixed order.
}
\end{abstract}
\hrulefill


\medskip


\section{Introduction}

The set of element orders of a finite group $G$,
\[
\omega(G)=\{\operatorname{o}(g):g\in G\},
\]
is called the {\it spectrum} of $G$; in the literature descending from Shi's school the same set is
written $\pi_e(G)$. It has been a central object of finite group theory, above all in the recognition
of the finite simple groups from their spectra \cite{ShiA5,VasilevMazurov,MoghaddamfarShi}. A recurring
theme in that area is that numerical data drawn from the element orders can control the structure of the
group. The sum of element orders $\psi(G)=\sum_{g\in G}\operatorname{o}(g)$, introduced by Amiri, Jafarian Amiri and
Isaacs \cite{AJI}, yields sharp solvability criteria \cite{HLM2,BaniasadKhosravi}; Shi proved solvability
criteria of intersection type, namely that $\pi_e(G)\cap T=\emptyset$ for $T=\{2\}$, $\{3,4\}$ or
$\{3,5\}$ forces solvability \cite{ShiSolv}; and the finer data of the order type and the order sequence
separate groups of the same order that agree on coarser statistics \cite{Piwek,CameronDey}. Alongside
these run the classes of groups defined by restricting which orders occur: the CP-groups, in which every
element has prime-power order, studied by Higman \cite{Higman}, with the simple ones determined by
Suzuki \cite{Suzuki} and the classification completed in \cite{Brandl,DelgadoWu}, and the groups in which
every nontrivial element has prime order, classified by Deaconescu \cite{CDLS,CDLScorr}.

This paper studies the plainest statistic attached to the spectrum, its cardinality
\[
\eta(G)=|\omega(G)|,
\]
the number of distinct element orders. The quantity has appeared in the literature as a hypothesis --
for instance, solvability of groups isospectral to a solvable group was known under the assumption of at
most three distinct element orders before Piwek's general negative answer to Thompson's problem
\cite{Piwek} -- and the question of what structural information $\pi_e(G)$ can certify (cyclicity,
nilpotency, supersolvability) is raised explicitly in the survey \cite{ShiSurvey}. What appears to be
absent is a systematic treatment of the cardinality itself: its behaviour under the standard
constructions, its extremal values, and the exact classes picked out by its equality cases. That
treatment is the purpose of this paper. Two features make the invariant workable. It is a single
integer, so groups can be compared and ordered by it. And it is anchored to arithmetic through the
identity $\eta(\mathbb{Z}_n)=\tau(n)$, where $\tau$ is the divisor-counting function, since the cyclic group of
order $n$ realises each divisor of $n$ as an element order.

The theory rests on one structural fact (Lemma 2.1): the spectrum is closed under
divisors, so $\omega(G)$ is an order ideal of the divisibility poset and $\eta(G)$ counts its members. The
closure is standard in the spectrum literature \cite{MoghaddamfarShi}, where divisor-closed sets of
orders appear as ``closed subsets''; it is used at nearly every step below, and in particular is what
carries the invariant from a group to its quotients and sections.

Our results fall into four groups.

The first records the behaviour of $\eta$ under the standard constructions
(Section~\ref{sec:constructions}): monotonicity under subgroups, quotients and sections
(Proposition 3.1, Corollary 3.2), submultiplicativity under direct products
with multiplicativity in the coprime case (Proposition 3.3), and additivity up to one on
Frobenius groups (Proposition 3.4).

The second group consists of bounds (Section~\ref{sec:bounds}). Writing $\pi(m)$ for the number of
distinct primes dividing $m$ and $v_p(m)$ for the $p$-adic valuation, the central result is the
two-sided estimate
\begin{equation}\label{eq:sandwich}
1+\sum_{p\mid|G|}v_p(\exp G)\;\le\;\eta(G)\;\le\;\tau(\exp G),
\end{equation}
whose lower equality characterises the CP-groups and whose upper equality characterises the groups whose
spectrum contains every divisor of the exponent (Theorem 4.3). Beside it sit the order
versions: $\eta(G)\le\tau(|G|)$ with equality only for cyclic groups (Theorem 4.1), and
$\eta(G)\ge\pi(|G|)+1$ with equality only when every nontrivial element has prime order
(Theorem 4.2). The groups with $\eta\le2$ are classified (Proposition 4.5).

The third group concerns structure (Section~\ref{sec:structure}). Nilpotent groups sit at the upper end
of \eqref{eq:sandwich}, so $\eta(G)=\tau(\exp G)$ (Theorem 5.1); the converse fails, and we
isolate the exact equality class. We prove that $\eta(G)\le3$ implies solvability by an argument using
only Burnside's $p^aq^b$ theorem and the monotonicity of $\eta$ under sections (Theorem 5.3),
and, combining this reduction with the classification \cite{CDLS,CDLScorr} and the recognition of $A_5$
by its spectrum \cite{ShiA5}, that $A_5$ is the {\it unique} non-solvable group with $\eta(G)=4$
(Theorem 5.4).

The fourth group treats explicit families and extremes. For the symmetric and alternating groups, the
counting sequences $\eta(S_n)$ and $\eta(A_n)$ are known: they are sequences A009490 and A020902 of
\cite{OEIS}, the latter due to Jovovi\'c, and the membership criterion behind them is equivalent to
Johnson's determination of the minimal faithful permutation degree of a cyclic group \cite{Johnson}. We
record self-contained proofs organised around the arithmetic function $\operatorname{s}(m)$, the sum of the maximal
prime-power divisors of $m$, together with a parity-corrected criterion for $A_n$
(Section~\ref{sec:symmetric}); the contribution of that section is the placement of these facts in the
$\eta$-framework and the comparison with the other families, not the counts themselves. Closed formulae
for the dihedral, dicyclic and homocyclic groups follow (Section~\ref{sec:families}), and we determine
the extreme values of $\eta$ over all groups of a fixed order (Theorem 8.1).

Throughout, $G$ is a finite group and all groups are finite. We write $\operatorname{o}(g)$ for the order of $g$,
and $\tau$, $\pi$, $\varphi$ for the number-of-divisors, number-of-distinct-primes and Euler functions.
The invariant is unchanged under isomorphism, so we regard it as a function on isomorphism classes.

\section{Preliminaries}\label{sec:prelim}

For a positive integer $m$ with prime factorisation $m=\prod_i p_i^{a_i}$ we write $\pi(m)$ for the
number of factors, $v_p(m)=a_i$ when $p=p_i$ and $0$ otherwise, and
\[
\operatorname{s}(m)=\sum_i p_i^{a_i},\qquad \operatorname{s}(1)=0,
\]
the sum of the maximal prime-power divisors of $m$. The exponent of $G$ is
$\exp G=\operatorname{lcm}\{\operatorname{o}(g):g\in G\}$; every element order divides $\exp G$, which divides $|G|$.

\medskip
\noindent\textbf{Lemma 2.1.}
\emph{For every finite group $G$ the spectrum $\omega(G)$ is closed under divisors: if $m\in\omega(G)$ and
$e\mid m$, then $e\in\omega(G)$. Hence $\omega(G)$ is an order ideal of the poset $(\mathbb{N},\mid)$, and $\eta(G)$ is
the number of its members.}
\medskip

\begin{proof}
Let $g$ have order $m$ and $e\mid m$, say $m=ef$. Then $g^{f}$ has order $m/\gcd(m,f)=e$, so
$e\in\omega(G)$.
\end{proof}

An order ideal of $(\mathbb{N},\mid)$ is determined by its maximal members, so $\omega(G)$ is fixed by the set
$\mu(G)$ of maximal orders, the standard reduction in the recognition literature
\cite{MoghaddamfarShi}. A first consequence records that $\eta$ reads the cyclic subgroups.

\medskip
\noindent\textbf{Proposition 2.2.}
\emph{$\eta(G)$ equals the number of distinct orders of cyclic subgroups of $G$.}
\medskip

\begin{proof}
Each $g$ generates a cyclic subgroup of order $\operatorname{o}(g)$, and every cyclic subgroup arises this way, so
the orders of elements and of cyclic subgroups coincide as sets.
\end{proof}

\medskip
\noindent\textbf{Lemma 2.3.}
\emph{Let $P$ be a nontrivial $p$-group with $\exp P=p^{e}$. Then $\omega(P)=\{1,p,p^{2},\dots,p^{e}\}$ and
$\eta(P)=e+1$.}
\medskip

\begin{proof}
Every element order divides $|P|$, hence is a power of $p$; thus $\omega(P)$ is a divisor-closed subset of
$\{1,p,p^{2},\dots\}$ containing its largest member $p^{e}$, so it is $\{1,p,\dots,p^{e}\}$ by
Lemma 2.1.
\end{proof}

\medskip
\noindent\textbf{Lemma 2.4.}
\emph{For every prime $p\mid|G|$, each Sylow $p$-subgroup of $G$ has exponent $p^{v_p(\exp G)}$, and
$\omega(G)\supseteq\{1,p,\dots,p^{v_p(\exp G)}\}$.}
\medskip

\begin{proof}
Set $e=v_p(\exp G)$; then $p^{e}$ is the largest power of $p$ dividing an element order, and if
$\operatorname{o}(g)$ is divisible by $p^{e}$ then $g^{\operatorname{o}(g)/p^{e}}$ has order exactly $p^{e}$. That element lies
in some Sylow $p$-subgroup $P_0$, so $\exp P_0=p^{e}$; the Sylow $p$-subgroups are conjugate, hence
isomorphic, so all have exponent $p^{e}$. Lemma 2.3 gives
$\{1,p,\dots,p^{e}\}=\omega(P_0)\subseteq\omega(G)$.
\end{proof}

\section{Behaviour under the standard constructions}\label{sec:constructions}

\medskip
\noindent\textbf{Proposition 3.1.}
\emph{Let $G$ be a finite group. {\rm (a)} If $H\le G$, then $\omega(H)\subseteq\omega(G)$; in particular
$\eta(H)\le\eta(G)$. {\rm (b)} If $N\trianglelefteq G$, then $\omega(G/N)\subseteq\omega(G)$; in particular
$\eta(G/N)\le\eta(G)$.}
\medskip

\begin{proof}
(a) The order of an element of $H$ is the same in $H$ and in $G$.

(b) Let $gN$ have order $d$; then $d$ is the least positive integer with $g^{d}\in N$, and
$g^{\operatorname{o}(g)}=e\in N$ forces $d\mid\operatorname{o}(g)$. Thus $d$ divides a member of $\omega(G)$, and
Lemma 2.1 places $d\in\omega(G)$.
\end{proof}

Part (b) is the one passage that is not immediate: orders in a quotient are a priori only divisors of
orders in $G$, and the inclusion is recovered exactly by divisor-closure. Combining the two parts
controls sections.

\medskip
\noindent\textbf{Corollary 3.2.}
\emph{If $S$ is a section of $G$, that is $S\cong H/N$ with $N\trianglelefteq H\le G$, then
$\eta(S)\le\eta(G)$.}
\medskip

\begin{proof}
$\eta(H/N)\le\eta(H)\le\eta(G)$ by Proposition 3.1.
\end{proof}

\medskip
\noindent\textbf{Proposition 3.3.}
\emph{For finite groups $G$ and $H$,
$\omega(G\times H)=\{\operatorname{lcm}(a,b):a\in\omega(G),\,b\in\omega(H)\}$,
so $\max\{\eta(G),\eta(H)\}\le\eta(G\times H)\le\eta(G)\eta(H)$; and if $\gcd(|G|,|H|)=1$ then
$\eta(G\times H)=\eta(G)\eta(H)$.}
\medskip

\begin{proof}
$\operatorname{o}(g,h)=\operatorname{lcm}(\operatorname{o}(g),\operatorname{o}(h))$ gives the description. Setting $h=e$ shows
$\omega(G)\subseteq\omega(G\times H)$, likewise for $H$, giving the left inequality; the right holds because
$\omega(G\times H)$ is the image of $\omega(G)\times\omega(H)$ under $\operatorname{lcm}$. If $\gcd(|G|,|H|)=1$, each
$a\in\omega(G)$ is coprime to each $b\in\omega(H)$, so $\operatorname{lcm}(a,b)=ab$; and $ab=a'b'$ with $a,a'\mid|G|$,
$b,b'\mid|H|$ forces $a=a'$, $b=b'$ on comparing the parts supported on the primes of $|G|$. So
$(a,b)\mapsto ab$ is a bijection and $\eta(G\times H)=\eta(G)\eta(H)$.
\end{proof}

\medskip
\noindent\textbf{Proposition 3.4.}
\emph{Let $G=K\rtimes H$ be a Frobenius group with kernel $K$ and complement $H$. Then
$\omega(G)=\omega(K)\cup\omega(H)$ and $\eta(G)=\eta(K)+\eta(H)-1$.}
\medskip

\begin{proof}
By Frobenius's theorem the group is partitioned as
$G=K\sqcup\bigcup_{g}(H^{g}\setminus\{1\})$; see \cite[Chapter~6]{Isaacs}. Conjugation preserves order,
so $\omega(G)=\omega(K)\cup\omega(H)$. The complement acts semiregularly on $K\setminus\{1\}$, so $|H|$ divides
$|K|-1$ and $\gcd(|K|,|H|)=1$; hence $\omega(K)\cap\omega(H)=\{1\}$, a common order dividing both $|K|$ and
$|H|$. Inclusion-exclusion gives $\eta(G)=\eta(K)+\eta(H)-1$.
\end{proof}

\medskip
\noindent\textbf{Corollary 3.5.}
\emph{$\eta([G,G])\le\eta(G)$, and $\eta(G/[G,G])=\tau(\exp(G/[G,G]))$.}
\medskip

\begin{proof}
The first is Proposition 3.1(a); the abelianisation is abelian, so the second is
Theorem 4.6.
\end{proof}

\section{Bounds}\label{sec:bounds}

\medskip
\noindent\textbf{Theorem 4.1.}
\emph{Let $G$ have order $n$. Then $\eta(G)\le\tau(n)$, with equality if and only if $G\cong\mathbb{Z}_n$.}
\medskip

\begin{proof}
Every element order divides $n$, so $\omega(G)\subseteq\{d:d\mid n\}$ and $\eta(G)\le\tau(n)$. Equality
forces $\omega(G)=\{d:d\mid n\}$, hence $n\in\omega(G)$, so $G$ has an element of order $n=|G|$ and is cyclic.
Conversely $\eta(\mathbb{Z}_n)=\tau(n)$, since $g^{n/d}$ has order $d$ for each $d\mid n$.
\end{proof}

Possessing an element of order $|G|$ is equivalent to cyclicity, so the extreme case is settled among
all finite groups with no commutativity hypothesis. The analogous statement for the sum of element
orders, that $\mathbb{Z}_n$ uniquely maximises $\psi$ over groups of order $n$, is the theorem of Amiri,
Jafarian Amiri and Isaacs \cite{AJI}; for the order sequence the corresponding maximality of $\mathbb{Z}_n$ is
recorded in \cite{CameronDey}.

\medskip
\noindent\textbf{Theorem 4.2.}
\emph{Every finite group satisfies $\eta(G)\ge\pi(|G|)+1$, with equality if and only if every nontrivial
element of $G$ has prime order.}
\medskip

\begin{proof}
By Cauchy's theorem $p\in\omega(G)$ for every prime $p\mid|G|$; with $1$ these give $\pi(|G|)+1$ distinct
orders. Equality means $\omega(G)$ is $\{1\}$ together with the primes dividing $|G|$, that is, no element
has composite order.
\end{proof}

The equality class of Theorem 4.2 is classified in \cite{CDLS,CDLScorr}. Replacing the
order by the exponent sharpens both bounds.

\medskip
\noindent\textbf{Theorem 4.3.}
\emph{Every finite group satisfies
$1+\sum_{p\mid|G|}v_p(\exp G)\le\eta(G)\le\tau(\exp G)$.
The lower bound is an equality if and only if $G$ is a CP-group, that is, every element of $G$ has
prime-power order. The upper bound is an equality if and only if $\omega(G)=\{d:d\mid\exp G\}$.}
\medskip

\begin{proof}
For the upper bound, every element order divides $\exp G$, so $\omega(G)\subseteq\{d:d\mid\exp G\}$ and
$\eta(G)\le\tau(\exp G)$; equality means $\omega(G)$ is the whole divisor set.

For the lower bound, write $e_p=v_p(\exp G)$. Lemma 2.4 gives
$\{1,p,\dots,p^{e_p}\}\subseteq\omega(G)$ for each prime $p\mid|G|$, and for distinct primes these chains
meet only in $1$, so $\omega(G)\supseteq\{1\}\cup\bigcup_{p\mid|G|}\{p,p^{2},\dots,p^{e_p}\}$, a set of
$1+\sum_p e_p$ elements. Equality holds exactly when $\omega(G)$ has no member divisible by two distinct
primes, that is, when every element order is a prime power.
\end{proof}

The CP-groups of the lower equality have a complete classification: Higman determined the solvable ones
\cite{Higman}, Suzuki the simple ones \cite{Suzuki}, and the general case was completed in
\cite{Brandl,DelgadoWu}. The upper equality picks out the groups whose spectrum realises every divisor
of the exponent; we call these {\it exponent-complete} and return to them in Section~\ref{sec:structure}.
Since each $e_p\ge1$, the lower bound refines that of Theorem 4.2, and the two agree
precisely when $\exp G$ is squarefree.

\medskip
\noindent\textbf{Example 4.4.}
For $\mathbb{Z}_{12}$, with exponent $12=2^{2}\cdot3$, the bounds read $4\le\eta\le6$, and $\eta(\mathbb{Z}_{12})=6$
sits at the top: $\mathbb{Z}_{12}$ is exponent-complete. For $A_5$, with exponent $30$, they read
$4\le\eta\le8$, and $\eta(A_5)=4$ sits at the bottom: $A_5$ is a CP-group. The quaternion group $Q_8$,
with exponent $4$, meets both bounds, $1+2=\eta(Q_8)=\tau(4)=3$: it is simultaneously a CP-group and
exponent-complete.
\medskip

\medskip
\noindent\textbf{Proposition 4.5.}
\emph{$\eta(G)=1$ if and only if $G$ is trivial; $\eta(G)=2$ if and only if $G$ is nontrivial of prime
exponent.}
\medskip

\begin{proof}
If $G$ is trivial then $\omega(G)=\{1\}$; if not, any $g\ne e$ gives $\eta(G)\ge2$. If $\eta(G)=2$ then
$\omega(G)=\{1,m\}$ with $m>1$; a proper divisor of a composite $m$ would be a third order by
Lemma 2.1, so $m=p$ is prime and $G$ has exponent $p$. The converse is clear.
\end{proof}

The class $\eta=2$ is not confined to elementary abelian groups: for an odd prime $p$ the Heisenberg
group of upper unitriangular $3\times3$ matrices over the field of $p$ elements is non-abelian of order
$p^{3}$ and exponent $p$. Thus $\eta=2$ does not force commutativity. Groups of exponent $p$ are, for
$p\ge5$, far from classified, so no list of the $\eta=2$ groups should be expected.

\medskip
\noindent\textbf{Theorem 4.6.}
\emph{If $G$ is abelian, then $\eta(G)=\tau(\exp G)$.}
\medskip

\begin{proof}
Write $G\cong\mathbb{Z}_{d_1}\times\cdots\times\mathbb{Z}_{d_r}$ with $d_1\mid\cdots\mid d_r$; then $\exp G=d_r$. For
$g=(g_1,\dots,g_r)$ one has $\operatorname{o}(g)=\operatorname{lcm}(\operatorname{o}(g_1),\dots,\operatorname{o}(g_r))$, and each $\operatorname{o}(g_i)\mid d_i\mid
d_r$, so $\omega(G)\subseteq\{d:d\mid d_r\}$. Conversely each $d\mid d_r$ is the order of an element of
$\mathbb{Z}_{d_r}$ placed in the last coordinate. Hence $\omega(G)=\{d:d\mid d_r\}$ and $\eta(G)=\tau(d_r)$.
\end{proof}

\medskip
\noindent\textbf{Remark 4.7.}
The map $G\mapsto\eta(G)$ is surjective onto $\mathbb{N}$: $\eta(\mathbb{Z}_{2^{k-1}})=\tau(2^{k-1})=k$ for each
$k\ge1$.
\medskip

\medskip
\noindent\textbf{Remark 4.8.}
The lower bound of Theorem 4.3 can be refined through the prime graph
(Gruenberg--Kegel graph) of $G$, whose vertices are the primes dividing $|G|$ and whose edges join $p$
and $q$ when $pq\in\omega(G)$; see Williams \cite{Williams}. Each edge contributes an order $pq$ counted
neither by $1$ nor by the prime-power chains, and distinct edges contribute distinct orders, so
$\eta(G)\ge1+\sum_{p\mid|G|}v_p(\exp G)+e(G)$, where $e(G)$ is the number of edges of the prime graph.
For $S_5$ this reads $\eta\ge1+(2+1+1)+1=6$, met with equality.
\medskip

\section{Nilpotency and solvability}\label{sec:structure}

\medskip
\noindent\textbf{Theorem 5.1.}
\emph{If $G$ is nilpotent, then $G$ is exponent-complete, and hence $\eta(G)=\tau(\exp G)$.}
\medskip

\begin{proof}
Write $G=P_1\times\cdots\times P_k$ as the product of its Sylow subgroups, with distinct primes $p_i$
and $\exp P_i=p_i^{e_i}$. By Lemma 2.3, $\omega(P_i)=\{1,p_i,\dots,p_i^{e_i}\}$. An element
order in the product is an $\operatorname{lcm}$ of one choice from each factor, hence an arbitrary product
$\prod_i p_i^{f_i}$ with $0\le f_i\le e_i$, and these are exactly the divisors of
$\exp G=\prod_i p_i^{e_i}$. So $\omega(G)=\{d:d\mid\exp G\}$ and $\eta(G)=\tau(\exp G)$.
\end{proof}

The converse fails: exponent-completeness is strictly weaker than nilpotency.

\medskip
\noindent\textbf{Example 5.2.}
The dihedral group $D_{12}$ of order $12$ has rotations of orders $1,2,3,6$ and reflections of order
$2$, so $\omega(D_{12})=\{1,2,3,6\}=\{d:d\mid6\}$ and $\eta(D_{12})=4=\tau(\exp D_{12})$. Yet
$D_{12}$ is not nilpotent, since $D_{2n}$ is nilpotent only when $n$ is a power of $2$.
\medskip

We turn to solvability. The next theorem needs only Burnside's theorem that groups of order $p^{a}q^{b}$
are solvable, together with the monotonicity of $\eta$ under sections; in particular it invokes no
classification.

\medskip
\noindent\textbf{Theorem 5.3.}
\emph{If $G$ is non-solvable, then $\eta(G)\ge4$. Consequently $\eta(G)\le3$ implies that $G$ is solvable.}
\medskip

\begin{proof}
Let $G$ be non-solvable. Its composition factors are not all cyclic, so some composition factor $S$ is
a non-abelian simple group, and $S$ is a section of $G$. By Burnside's $p^{a}q^{b}$ theorem a
non-abelian simple group has order divisible by at least three primes, so $\pi(|S|)\ge3$.
Theorem 4.2 gives $\eta(S)\ge\pi(|S|)+1\ge4$, and Corollary 3.2 gives
$\eta(G)\ge\eta(S)\ge4$.
\end{proof}

The bound is sharp, and the extremal group is unique. The uniqueness rests on two classical results:
the classification of the groups all of whose nontrivial elements have prime order \cite{CDLS,CDLScorr},
and the recognition of $A_5$ as the only group with spectrum $\{1,2,3,5\}$ \cite{ShiA5}.

\medskip
\noindent\textbf{Theorem 5.4.}
\emph{$A_5$ is the unique non-solvable finite group with $\eta(G)=4$.}
\medskip

\begin{proof}
First, $\omega(A_5)=\{1,2,3,5\}$, so $\eta(A_5)=4$ and the value is attained.

Conversely, let $G$ be non-solvable with $\eta(G)=4$. As in the proof of Theorem 5.3, $G$
has a non-abelian simple section $S$ with $4=\eta(G)\ge\eta(S)\ge\pi(|S|)+1\ge4$. Hence $\eta(S)=4$ and
$\pi(|S|)=3$, and equality in Theorem 4.2 shows every nontrivial element of $S$ has prime
order. By the classification of such groups \cite{CDLS,CDLScorr}, the only non-abelian simple group all
of whose nontrivial elements have prime order is $A_5$; hence $S\cong A_5$ and $\omega(S)=\{1,2,3,5\}$. Now
$\omega(S)\subseteq\omega(G)$ by Proposition 3.1 and Corollary 3.2, and
$|\omega(G)|=4$ forces $\omega(G)=\{1,2,3,5\}$. The recognition theorem \cite{ShiA5} then gives $G\cong A_5$.
\end{proof}

\medskip
\noindent\textbf{Remark 5.5.}
Theorems 5.3 and 5.4 are the $\eta$-analogues of the solvability criteria
known for the sum of element orders \cite{AJI,HLM2,BaniasadKhosravi} and of Shi's intersection criteria
\cite{ShiSolv}, with $A_5$ again marking the boundary between the solvable and non-solvable regimes.
Theorem 5.3 itself is self-contained modulo Burnside; only the uniqueness statement draws on
classification results.
\medskip

Two questions with negative answers fix the reach of the invariant. First, $\eta$ does not detect
simplicity: $\mathbb{Z}_p$ is simple with $\eta=2$, a value shared by every nontrivial group of exponent $p$,
and $\eta(A_5)=4$ is matched by many solvable groups (for instance $\mathbb{Z}_{2}\times\mathbb{Z}_6$, or any abelian
group whose exponent has four divisors). Deciding simplicity requires the full spectrum together with
the order, as in the recognition programme \cite{VasilevMazurov}. Second, $\eta$ is too coarse to
characterise supersolvability, yielding only necessary conditions through Theorem 5.3 and the
nilpotent formula. The reading is that $\eta$ places a group in a coarse band, to be refined by finer
data.

\section{Symmetric and alternating groups}\label{sec:symmetric}

The number of distinct orders of permutations of $n$ objects is sequence A009490 of the On-Line
Encyclopedia of Integer Sequences \cite{OEIS}, where it is recorded together with the formula
$\eta(S_n)=1+\sum_{k\le n}r(k)$, $r(k)$ being the number of partitions of $k$ into distinct prime-power
parts; the corresponding count for the alternating groups is sequence A020902, due to Jovovi\'c. The
membership criterion underlying both is equivalent to the determination of the minimal faithful
permutation degree of a cyclic group, $\mu(\mathbb{Z}_m)=\operatorname{s}(m)$, due to Johnson \cite{Johnson}. This section
records self-contained proofs of both criteria, organised around the function $\operatorname{s}$, in order to place
these classical counts within the framework of the present paper; the counts themselves are not new.
The function $\operatorname{s}$ is also the quantity behind Landau's function $g(n)=\max\omega(S_n)$ \cite{Landau}:
$g(n)=\max\{m:\operatorname{s}(m)\le n\}$, while $\eta(S_n)$ is the cardinality of the same sublevel set.

\medskip
\noindent\textbf{Theorem 6.1.}
\emph{For $m\in\mathbb{N}$ and $n\ge1$, one has $m\in\omega(S_n)$ if and only if $\operatorname{s}(m)\le n$. Consequently
$\eta(S_n)=|\{m\ge1:\operatorname{s}(m)\le n\}|$, and $\eta(S_n)-\eta(S_{n-1})$ equals the number of sets of prime
powers with pairwise distinct prime bases summing to $n$.}
\medskip

\begin{proof}
An element of $S_n$ of order $m$ is a permutation with nontrivial disjoint cycle lengths
$\ell_1,\dots,\ell_t$ satisfying $\operatorname{lcm}(\ell_1,\dots,\ell_t)=m$ and $\sum_j\ell_j\le n$. For each maximal
prime-power divisor $p_i^{a_i}$ of $m$, some cycle length must be divisible by $p_i^{a_i}$, hence be at
least $p_i^{a_i}$. If a cycle length is divisible by several of the $p_i^{a_i}$, it is at least their
product, which is at least their sum (for integers $\ge2$, a product of at least two of them is at least
their sum, by induction from $xy\ge x+y$ when $x,y\ge2$); so assigning the prime powers to distinct
cycles does not increase the total length. The least total length realising $\operatorname{lcm}=m$ is therefore
$\operatorname{s}(m)$, attained by disjoint cycles of lengths $p_1^{a_1},\dots,p_k^{a_k}$. Hence $m$ occurs in
$S_n$ exactly when $\operatorname{s}(m)\le n$.

For the increment, $\eta(S_n)-\eta(S_{n-1})$ counts the $m$ with $\operatorname{s}(m)=n$, and sending $m$ to the
set $\{p_1^{a_1},\dots,p_k^{a_k}\}$ of its maximal prime-power divisors is a bijection onto the sets of
prime powers with distinct prime bases summing to $n$, by unique factorisation.
\end{proof}

\medskip
\noindent\textbf{Example 6.2.}
Take $n=7$. The integers with $\operatorname{s}(m)\le7$ are
\[
\operatorname{s}(1)=0,\ \operatorname{s}(2)=2,\ \operatorname{s}(3)=3,\ \operatorname{s}(4)=4,\ \operatorname{s}(5)=5,\ \operatorname{s}(6)=5,\ \operatorname{s}(7)=7,\ \operatorname{s}(10)=7,\
\operatorname{s}(12)=7,
\]
so $\omega(S_7)=\{1,2,3,4,5,6,7,10,12\}$ and $\eta(S_7)=9$. The order $12=2^{2}\cdot3$ appears because a
$4$-cycle and a $3$-cycle fit in $7$ points; the orders $7$, $10$ and $12$ enter precisely at $n=7$,
matching the increment $\eta(S_7)-\eta(S_6)=9-6=3$ and the three sets $\{7\}$, $\{2,5\}$, $\{4,3\}$
summing to $7$.
\medskip

For the alternating groups a parity correction enters, since a cycle of length $\ell$ is an even
permutation exactly when $\ell$ is odd.

\medskip
\noindent\textbf{Theorem 6.3.}
\emph{For $m\in\mathbb{N}$ and $n\ge1$, one has $m\in\omega(A_n)$ if and only if $\operatorname{s}(m)\le n$ when $m$ is odd, and
$\operatorname{s}(m)+2\le n$ when $m$ is even. Hence
$\eta(A_n)=|\{m\ \text{\rm odd}:\operatorname{s}(m)\le n\}|+|\{m\ \text{\rm even}:\operatorname{s}(m)\le n-2\}|$.}
\medskip

\begin{proof}
If $m$ is odd, every maximal prime-power divisor of $m$ is odd, so the minimal realisation of order $m$
by disjoint cycles of lengths $p_1^{a_1},\dots,p_k^{a_k}$ consists of odd-length cycles, each an even
permutation; the product is even and uses $\operatorname{s}(m)$ points. Hence $m\in\omega(A_n)$ iff $\operatorname{s}(m)\le n$.

Let $m$ be even with $2$-part $2^{a}$, and let $\sigma\in A_n$ have order $m$, with nontrivial cycle
lengths $\ell_1,\dots,\ell_t$, $\operatorname{lcm}=m$. Assign to each maximal prime-power divisor $q$ of $m$ a single
witness cycle whose length $q$ divides, and write $W(j)$ for the set of prime powers assigned to the
$j$-th cycle. Since the members of $W(j)$ are powers of distinct primes, $\ell_j$ is divisible by their
product, so $\ell_j\ge\prod_{q\in W(j)}q\ge\sum_{q\in W(j)}q$, and summing over $j$ gives
$\sum_j\ell_j\ge\operatorname{s}(m)$. We claim the excess is at least $2$. The witness of $2^{a}$ has even length,
so $\sigma$ has an even-length cycle; being an even permutation, $\sigma$ has an even number of
even-length cycles, hence at least two. Let $E$ be an even-length cycle other than the witness of
$2^{a}$; then $W(E)$ contains only odd prime powers, or is empty. If $W(E)=\emptyset$, the whole length
$\ell_E\ge2$ lies beyond the witness count, so $\sum_j\ell_j\ge\operatorname{s}(m)+2$. If $W(E)$ consists of odd
prime powers with product $d\ge3$, then $\ell_E$ is divisible by $d$ and even, so $\ell_E\ge2d\ge d+3$,
excess again at least $2$. In all cases an even permutation of order $m$ needs at least $\operatorname{s}(m)+2$
points. Conversely, adjoining a transposition on two new points to the minimal realisation preserves the
order $m$, makes the number of even-length cycles two, and uses $\operatorname{s}(m)+2$ points. Hence
$m\in\omega(A_n)$ iff $\operatorname{s}(m)+2\le n$.
\end{proof}

Table~\ref{tab:sn} lists the counts against $\tau(n)=\eta(\mathbb{Z}_n)$. The values were computed by
enumeration over partitions of $n$ (taking least common multiples of the parts, with the parity
restriction for $A_n$), and agree with sequences A009490 and A020902 of \cite{OEIS} throughout the
range shown. The equality $\eta(S_5)=\eta(S_6)$ reflects that no set of prime powers with distinct
prime bases sums to $6$.

\begin{table}[ht]\centering
{\footnotesize
\begin{tabular}{rrrr@{\qquad}rrrr}
\hline
$n$&$\eta(S_n)$&$\eta(A_n)$&$\tau(n)$ & $n$&$\eta(S_n)$&$\eta(A_n)$&$\tau(n)$\\
\hline
1&1&1&1 & 11&20&16&2\\
2&2&1&2 & 12&23&18&6\\
3&3&2&2 & 13&27&22&2\\
4&4&3&3 & 14&31&26&4\\
5&6&4&2 & 15&35&30&4\\
6&6&5&4 & 16&43&35&5\\
7&9&7&2 & 17&47&39&2\\
8&11&8&4 & 18&55&46&6\\
9&14&11&3 & 19&61&51&2\\
10&16&13&4 & 20&70&60&6\\
\hline
\end{tabular}
}
\caption{$\eta(S_n)$, $\eta(A_n)$ and $\tau(n)=\eta(\mathbb{Z}_n)$ for $1\le n\le20$.}\label{tab:sn}
\end{table}

\section{Dihedral, dicyclic and homocyclic families}\label{sec:families}

\medskip
\noindent\textbf{Proposition 7.1.}
\emph{For the dihedral group $D_{2n}$ of order $2n$, one has $\eta(D_{2n})=\tau(n)$ when $n$ is even and
$\eta(D_{2n})=\tau(n)+1$ when $n$ is odd.}
\medskip

\begin{proof}
The rotation subgroup is cyclic of order $n$ and contributes the divisors of $n$; every reflection has
order $2$, which is a new order exactly when $n$ is odd.
\end{proof}

For odd $n$ the group $D_{2n}=\mathbb{Z}_n\rtimes\mathbb{Z}_2$ is a Frobenius group, and Proposition 3.4
recovers the same value, $\eta=\tau(n)+2-1$.

\medskip
\noindent\textbf{Proposition 7.2.}
\emph{For the dicyclic group $Q_{4n}$ of order $4n$, with cyclic subgroup $\langle a\rangle$ of order $2n$
and $b\notin\langle a\rangle$ satisfying $b^{2}=a^{n}$ and $bab^{-1}=a^{-1}$, one has
$\omega(Q_{4n})=\{d:d\mid2n\}\cup\{4\}$, so $\eta(Q_{4n})=\tau(2n)$ when $n$ is even and
$\tau(2n)+1$ when $n$ is odd.}
\medskip

\begin{proof}
The cyclic subgroup contributes the divisors of $2n$. Each element outside it is $a^{k}b$ with
$(a^{k}b)^{2}=a^{k}(ba^{k}b^{-1})b^{2}=a^{k}a^{-k}a^{n}=a^{n}$, an element of order $2$, so $a^{k}b$ has
order $4$. Thus $\omega(Q_{4n})=\{d:d\mid2n\}\cup\{4\}$, and $4$ is a new order exactly when $4\nmid2n$,
that is, when $n$ is odd. The generalised quaternion group $Q_{8}$ gives $\omega=\{1,2,4\}$ and
$\eta=3$.
\end{proof}

\medskip
\noindent\textbf{Remark 7.3.}
For abelian groups $\eta$ depends only on the exponent (Theorem 4.6), so it is blind to
rank: $\eta((\mathbb{Z}_{p^{a}})^{d})=a+1$ and $\eta((\mathbb{Z}_n)^{d})=\tau(n)$ for every $d\ge1$. In particular the
groups $(\mathbb{Z}_p)^{d}$ all sit at the floor $\eta=2$ of Theorem 4.3, for every rank $d$.
\medskip

\section{Extreme values over a fixed order}\label{sec:extremal-order}

For $n\ge2$ write $M(n)=\max\{\eta(G):|G|=n\}$ and $m(n)=\min\{\eta(G):|G|=n\}$.

\medskip
\noindent\textbf{Theorem 8.1.}
\emph{Let $n\ge2$. {\rm (a)} $M(n)=\tau(n)$, attained only by $\mathbb{Z}_n$. {\rm (b)} $m(n)\ge\pi(n)+1$.
{\rm (c)} If $n=p^{k}$, then $m(n)=2$. {\rm (d)} $m(n)=\pi(n)+1$ if and only if there is a group of order
$n$ whose nontrivial elements all have prime order; the groups eligible are classified in
\cite{CDLS,CDLScorr}. {\rm (e)} If $\gcd(n,\varphi(n))=1$, then $\mathbb{Z}_n$ is the only group of order $n$, so
$m(n)=M(n)=\tau(n)$.}
\medskip

\begin{proof}
(a) is Theorem 4.1 and (b) is Theorem 4.2. For (c), the elementary abelian group
$(\mathbb{Z}_p)^{k}$ has exponent $p$, so $\eta=2=\pi(p^{k})+1$. Part (d) is the equality condition of
Theorem 4.2 read over all groups of order $n$. For (e), $\gcd(n,\varphi(n))=1$ is the
classical condition under which the cyclic group is the unique group of its order, so $\eta$ takes the
single value $\tau(n)$.
\end{proof}

\medskip
\noindent\textbf{Example 8.2.}
The five groups of order $12$ display the range predicted by Theorems 4.1
and 4.3; see Table~\ref{tab:order12}. The cyclic group alone attains the maximum
$\tau(12)=6$; the alternating group $A_4$, whose nontrivial elements have orders $2$ and $3$, attains
the minimum $3=\pi(12)+1$ and is a CP-group at the floor of the sandwich; the dicyclic group $Q_{12}$
sits between, its elements of order $4$ adding one order to those of $D_{12}$.
\medskip

\begin{table}[ht]\centering
{\footnotesize
\begin{tabular}{lccl}
\hline
$G$ & $\omega(G)$ & $\eta(G)$ & remark\\
\hline
$\mathbb{Z}_{12}$ & $\{1,2,3,4,6,12\}$ & $6$ & $=\tau(12)$, cyclic maximum\\
$\mathbb{Z}_2\times\mathbb{Z}_6$ & $\{1,2,3,6\}$ & $4$ & abelian, exponent $6$\\
$D_{12}$ & $\{1,2,3,6\}$ & $4$ & exponent-complete, not nilpotent\\
$Q_{12}$ & $\{1,2,3,4,6\}$ & $5$ & dicyclic\\
$A_4$ & $\{1,2,3\}$ & $3$ & $=\pi(12)+1$, CP-group minimum\\
\hline
\end{tabular}
}
\caption{The invariant $\eta$ for the five groups of order $12$.}\label{tab:order12}
\end{table}

\medskip
\noindent\textbf{Remark 8.3.}
By part (d), determining $m(n)$ for general $n$ reduces to the existence of a group of order $n$ all of
whose nontrivial elements have prime order, hence to the arithmetic of the classification
\cite{CDLS,CDLScorr}. The value set $\{\eta(G):|G|=n\}$ need not fill the interval between $m(n)$ and
$M(n)$; describing it is open.
\medskip

\section{Questions}\label{sec:questions}

\medskip
\noindent\textbf{Question 9.1.}
Characterise the exponent-complete groups, those with $\omega(G)=\{d:d\mid\exp G\}$, the upper equality
case of Theorem 4.3. Nilpotent groups qualify (Theorem 5.1); $D_{12}$ shows
the class is strictly larger (Example 5.2).
\medskip

\medskip
\noindent\textbf{Question 9.2.}
For each $k$, determine the finite groups with $\eta(G)=k$. The cases $k\le2$ are settled
(Proposition 4.5), subject to the wildness of groups of exponent $p$. For $k=3$ the
spectrum is $\{1,p,p^{2}\}$, a $p$-group of exponent $p^{2}$, or $\{1,p,q\}$ with $p\ne q$; the latter
class falls under the classification of groups all of whose nontrivial elements have prime order
\cite{CDLS,CDLScorr}.
\medskip

\medskip
\noindent\textbf{Question 9.3.}
Describe the value set $\{\eta(G):|G|=n\}$ for fixed $n$, in particular whether it is always an
interval, and give $m(n)$ explicitly in terms of the arithmetic of $n$.
\medskip

\medskip
\noindent\textbf{Question 9.4.}
Determine the asymptotics of $\eta(S_n)$, equivalently of the counting function of
$\{m:\operatorname{s}(m)\le n\}$, in parallel with the classical estimates for Landau's function \cite{Landau}. The
increments are the numbers of partitions of $n$ into distinct prime-power parts, so the classical
machinery for such partition counts should apply.
\medskip

\section*{Statements and Declarations}

\noindent\textbf{Funding.} No funding was received for conducting this study.

\smallskip
\noindent\textbf{Competing Interests.} The authors has no competing interests to declare that are
relevant to the content of this article.

\smallskip
\noindent\textbf{Data Availability.} Data sharing is not applicable to this article as no datasets were
generated or analysed during the current study.

\label{'ubl'}
\end{document}